\documentclass[10pt]{article}
\usepackage{amsmath,amssymb,amsthm,amscd}
\numberwithin{equation}{section}

\newtheorem{prop}{Proposition}[section]
\newtheorem{theorem}[prop]{Theorem}

\newtheorem{remark}[prop]{Remark}

\newtheorem{acknowledgment}[prop]{Acknowledgment}

\def\<{\langle}
\def\>{\rangle}
\def\({\left(}
\def\){\right)}

\def\p{\partial}
\def\Ric{{\rm Ric}}

\begin{document}

\title{Scalar Curvature Behavior for Finite Time Singularity 
of K\"ahler-Ricci Flow}
\author{Zhou Zhang} 
\maketitle

\begin{abstract}

In this short paper, we show that K\"ahler-Ricci flows over 
closed manifolds would have scalar curvature blown-up 
for finite time singularity. Certain control of the blowing-up 
is achieved with some mild assumption. 

\end{abstract}

\section{Introduction}

In this short note, we consider the following K\"ahler-Ricci 
flow 
\begin{equation}
\label{krf}
\frac{\p\widetilde\omega_t}{\p t}=-{\rm Ric}(\widetilde\omega
_t)-\widetilde\omega_t, ~~~~\tilde\omega_0=\omega_0.
\end{equation}
over a closed K\"ahler manifold $X$ where $\omega_0$ is 
any K\"ahler metric on $X$. The short time existence of the 
solution is known from either R. Hamilton's general existence 
result on Ricci flow in \cite{ham} or the fact that K\"ahler-Ricci 
flow is indeed parabolic. \\ 

By the optimal existence result on K\"ahler-Ricci flow as in 
\cite{cascini-lanave} or \cite{t-znote}, we know the classic 
solution of (\ref{krf}) exists exactly as long as the cohomology 
class $[\widetilde\omega_t]$ from formal computation remains 
to be K\"ahler. The actually meaning will be explained later. 

It then comes down to the study of the behavior of the metric 
solution when approaching the flow singularity. In this work, 
we focus on the study of the case when the flow singularity 
happens at some finite time. Let's state the main theorems 
below.  

\begin{theorem}
\label{th:R-blows-up}
K\"ahler-Ricci flow (\ref{krf}) either exists for all time, or 
the scalar curvature blows up (from above) at some finite 
time, i. e.  
$$sup_{X\times[0, T)}|R(\widetilde\omega_t)|=+\infty$$   
where $T$ is the finite singular time. 
\end{theorem}

The blow-up would be from above in sight of the classic 
result on the lower bound of scalar curvature. Let's point 
out that the statement of this theorem also holds for other 
usual versions K\"ahler-Ricci flows by simple rescaling 
consideration. 

\begin{remark}

There have been some fundamental results regarding 
the finite time blowing-up of Ricci flow. In fact, it's known 
that curvature operator blows up by R. Hamilton's work 
\cite{ham4} and Ricci curvature blows up by N. Sesum's 
work \cite{sesum}. 

\end{remark} 
 
In the case of finite time singularity, suppose we have a 
holomorphic map 
$$F: X\longrightarrow Y$$ 
where $Y$ is an analytic variety smooth near the image 
$F(X)$ and there is a K\"ahler metric, $\omega_{_{M}}$ 
in a neighborhood of $F(X)$ such that $[\widetilde\omega
_t]=[F^*\omega_{_{M}}]$. Then we have the following 
control of the blowing-up of scalar curvature. 

\begin{theorem}
\label{th:R-control}
In the above setting, for the flow (\ref{krf}), 
$$R(\widetilde\omega_t)\leqslant\frac{C}{(T-t)^2}$$ 
where $C$ is a positive constant. 
\end{theorem} 

The motivation of the setting of this theorem is the 
semi-ampleness of the cohomology limit at singular 
time. It is of quite some interest in algebraic geometry 
as explained in \cite{t-znote}, for example. 

\begin{acknowledgment}

The author would like to thank R. Lazarsfeld for pointing 
out the result by J. Demailly and M. Paun in \cite
{demailly-paun} which is absolutely crucial to conclude  
Theorem \ref{th:R-blows-up} for general closed (K\"ahler) 
manifolds. The discussion with J. Song is also valuable 
for this result. The comparison with the result by G. 
Perelman as mentioned in the last section is suggested 
by J. Lott. The author can not thank his advisor, G. Tian, 
enough for introducing him to this interesting topic and 
constant encouragement along the way.  

\end{acknowledgment}

\section{Proof of Theorem \ref{th:R-blows-up}}

The proof is by contradiction. We assume {\bf the 
scalar curvature is uniformly bounded along the 
flow with finite time singularity}. One then makes 
use of J. Song and G. Tian's parabolic Schwarz 
Lemma (as in \cite{song-tian}) and some basic 
computations on this K\"ahler-Ricci flow to get 
some uniform control of the flow metric. The 
contradiction then comes from the general result 
on the existence of K\"ahler-Ricci flow and the 
numerical characterization of K\"ahler cone on 
closed K\"ahler manifolds by J. Demailly and M. 
Paun. The rest of this section contains the detailed 
argument. \\

As usual when dealing with K\"ahler-Ricci flow, we 
need the scalar version of (\ref{krf}) described below 
as in \cite{t-znote}. 

Let $\omega_t=\omega_\infty+e^{-t}(\omega_0-
\omega_\infty)$ where $[\omega_\infty]=K_X$ and 
$\omega_\infty=-{\rm Ric}(\Omega)$ for some smooth 
volume form $\Omega$ over $X$. Then set $\tilde
\omega_t=\omega_t+\sqrt{-1}\p\bar\p u$ and one 
has the following parabolic evolution equation for 
$u$,
\begin{equation}
\label{skrf}
\frac{\p u}{\p t}={\rm log}\frac{(\omega_t+\sqrt{-1}\p\bar
\p u)^n}{\Omega}-u, ~~~~u(\cdot, 0)=0.
\end{equation}

Now we state the following optimal existence result of 
K\"ahler-Ricci flow (as in \cite{cascini-lanave} and \cite
{t-znote}) mentioned in Introduction.  
\begin{prop}
\label{existence}
(\ref{krf}) (or (\ref{skrf}) equivalently) exists as long as 
$[\omega_t]$ remains K\"ahler, i. e. the solution is for 
the time interval $[0, T)$ where 
$$T=sup\{t|\,[\omega_t]~ \text{is K\"ahler}.\}.$$
\end{prop}

The finite time singularity means $[\omega_T]$ is on 
the boundary of the (open) K\"ahler cone, and thus no 
longer K\"ahler. Clearly it's "numerically effective" using 
the natural generalization of the notion from algebraic 
geometry. 

From now on, we consider the flow existing only for some 
finite interval $[0, T)$. As usual, the $C$'s below might 
stand for different positive constants. In case that the 
situation is more subtle, lower indices are used to tell 
them apart. The argument is organized into three steps. \\ 

\noindent $\bullet$ {\bf Volume Form (Lower) Bound} 

\vspace{0.1in}

With the bounded scalar curvature assumption, we 
can easily derive the uniform control on the volume 
form along the flow, using the following evolution of 
volume form,   
\begin{equation*}
\begin{split}
\frac{\p\tilde\omega^n_t}{\p t}
&= n\frac{\p\widetilde\omega_t}{\p t}\wedge{\tilde
\omega^{n-1}_t} \\
&= n(-{\rm Ric}(\widetilde\omega_t)-\widetilde
\omega_t)\wedge\tilde\omega^{n-1}_t \\
&= (-R-n)\tilde\omega^n_t. 
\end{split}
\end{equation*}

This gives $|\frac{\p u}{\p t}+u|\leqslant C$ as $\widetilde\omega^n_t=e^{\frac{\p u}{\p t}+u}\Omega$. 

\begin{remark}

Instead of the assumption on scalar curvature, one can 
also directly assume positive lower bound for the volume 
form or equivalently,  $\frac{\p u}{\p t}\geqslant -C$ since 
we are considering the finite time singularity case. This 
simple observation actually brings up a very intuitive 
understanding of Theorem \ref{th:R-blows-up}: the flow 
(\ref{skrf}) is stopped at some finite time because the 
term in $\log$ is tending to $0$, i. e. no uniform lower 
bound. 

\end{remark}

\vspace{0.1in}

\noindent $\bullet$ {\bf Metric Estimate}

\vspace{0.1in}

We begin with the inequality from parabolic Schwarz 
Lemma. In this note, the Laplacian $\Delta$ without 
lower index, is always with respect to the changing 
metric along the flow, $\widetilde\omega_t$. 

Let $\phi=\<\widetilde\omega_t, \omega_0\>>0$. Using 
computation for (\ref{krf}) in \cite{song-tian}, one has 
\begin{equation}
\label{pschwarz}
(\frac{\p}{\p t}-\Delta){\rm log}\phi\leqslant C_1\phi+1,
\end{equation}
where $C_1$ is a positive constant depending on the 
bisectional curvature of $\omega_0$. It's quite irrelevant 
that $\omega_0$ is the initial metric for the K\"ahler-Ricci 
flow. In fact, it doesn't have to be a metric over $X$ which 
is an interesting part of this lemma as indicated in \cite
{song-tian}. This is useful for the proof of Theorem 
\ref{th:R-control}. \\

Applying Maximum Principle to (\ref{skrf}) gives $u
\leqslant C$. Take $t$-derivative to get 
$$\frac{\p}{\p t}\left(\frac{\p u}{\p t}\right)=\Delta\left(
\frac{\p u}{\p t}\right)-e^{-t}\<\widetilde\omega_t, \omega
_0-\omega_\infty\>-\frac{\p u}{\p t}.$$
It can be reformulated into the following two equations, 
$$\frac{\p}{\p t}\left(e^t\frac{\p u}{\p t}\right)=\Delta\left(
e^t\frac{\p u}{\p t}\right)-\<\widetilde\omega_t, \omega_0-
\omega_\infty\>,$$
$$\frac{\p}{\p t}\left(\frac{\p u}{\p t}+u\right)=\Delta\left(
\frac{\p u}{\p t}+u\right)-n+\<\widetilde\omega_t, \omega_
\infty\>.$$
Their difference gives
\begin{equation}
\label{u-decreasing}
\frac{\p}{\p t}\left((e^t-1)\frac{\p u}{\p t}-u\right)=\Delta
\left((e^t-1)\frac{\p u}{\p t}-u\right)+n-\<\widetilde\omega_t, 
\omega_0\>.
\end{equation}
By Maximum Principle, this gives
$$(e^t-1)\frac{\p u}{\p t}-u-nt\leqslant C,$$
which together with the upper bound of $u$ and local 
bound for $\frac{\p u}{\p t}$ near $t=0$ would provide 
$$\frac{\p u}{\p t}\leqslant C.$$
The upper bounds on $\frac{\p u}{\p t}$ and $u$ together 
with $|\frac{\p u}{\p t}+u|\leqslant C$ from volume control 
give the uniform (lower) bounds on $\frac{\p u}{\p t}$ and 
$u$. \\

Multiply (\ref{u-decreasing}) by a large enough constant 
$C_2>C_1+1$ and combining it with (\ref{pschwarz}), one 
arrives at 
\begin{equation}
\label{main}
\begin{split}
(\frac{\p}{\p t}-\Delta)\left({\rm log}\phi+
(e^t-1)\frac{\p u}{\p t}-u\right)
&\leqslant nC_2+1-(C_2-C_1)\phi \\
&\leqslant C_3-\phi.
\end{split}
\end{equation}
Apply Maximum Principle for ${\rm log}\phi+(e^t-1)\frac
{\p u}{\p t}-u$. Considering the place where it achieves 
maximum value, one has 
$$\phi\leqslant C,$$
and so 
$${\rm log}\phi+(e^t-1)\frac{\p u}{\p t}-u\leqslant C.$$ 

Hence we conclude $\<\widetilde\omega_t, \omega_0\>
\leqslant C$ using the bound on $\frac{\p u}{\p t}$ 
and $u$. This trace bound, together with $\tilde
\omega^n_t\leqslant C\omega^n_0$, provide the 
uniform bound of $\widetilde\omega_t$ as metric, i. e. 
$C^{-1}\omega_0\leqslant \widetilde\omega_t \leqslant 
C\omega_0$. \\

\noindent $\bullet$ {\bf Contradiction}

\vspace{0.1in}

The metric (lower) bound makes sure that for any fixed 
analytic variety in $X$, the integral of the proper power 
of $\widetilde\omega_t$ is bounded away from $0$, and so 
the limiting class $[\omega_T]$ would have positive 
intersection with any analytic variety by taking the 
cohomology limit. Then by Theorem 4.1 in \cite
{demailly-paun}, we conclude that $[\omega_T]$ is actually 
K\"ahler which contradicts with the assumption of finite 
time singularity at $T$ in sight of Proposition 
\ref{existence}. \\

Hence we have finished the proof of Theorem \ref
{th:R-blows-up}. 

\begin{remark}

In sight of this numerical characterization of K\"ahler 
cone for any general closed K\"ahler manifold by J. 
Demailly and M. Paun, the blowing-up of curvature 
operator or Ricci curvature in closed K\"ahler case 
is fairly obvious. The situation of scalar curvature is 
the first non-trivial statement.   

\end{remark}

\section{Proof of Theorem \ref{th:R-control}}

Now we derive certain control of the blowing-up of 
scalar curvature by mainly following the argument 
in \cite{b-curvature}. The argument is also organized 
in three steps. \\ 

\noindent $\bullet$ {\bf $0$-th Order Estimates} 

\vspace{0.1in}

$u\leqslant C$ is directly from (\ref{skrf}). $t$-derivative 
of (\ref{skrf}) is  
$$\frac{\partial}{\partial t}(\frac{\partial u}{\partial t})=\Delta
(\frac{\partial u}{\partial t})-e^{-t}\<\widetilde\omega_t, \omega_0
-\omega_\infty\>-\frac{\partial u}{\partial t},$$
which has the following variations,
$$\frac{\partial}{\partial t}(e^t\frac{\partial u}{\partial t})=\Delta
(e^t\frac{\partial u}{\partial t})-\<\widetilde\omega_t, \omega_0-
\omega_\infty\>,$$ 
$$\frac{\partial}{\partial t}(\frac{\partial u}{\partial t}+u)=\Delta
(\frac{\partial u}{\partial t}+u)-n+\<\widetilde\omega_t, \omega_
\infty\>.$$
A proper linear combination of these equations provides 
the following "finite time version" of the second equation,
$$\frac{\partial}{\partial t}\((1-e^{t-T})\frac{\partial u}{\partial t}
+u\)=\Delta\((1-e^{t-T})\frac{\partial u}{\partial t}+u\)-n+\<\tilde
\omega_t, \omega_T\>.$$
The difference of the original two equations gives
$$\frac{\partial}{\partial t}\((1-e^t)\frac{\partial u}{\partial t}
+u\)=\Delta\((1-e^t)\frac{\partial u}{\partial t}+u\)-n+\<\tilde
\omega_t, \omega_0\>,$$
which implies the "essential decreasing" of metric potential 
along the flow, i. e.  
$$\frac{\partial u}{\partial t}\leqslant\frac{nt+C}{e^t-1}.$$

Notice that this estimate only depends on the initial value 
of $u$ and its upper bound along the flow. It is uniform 
away from the initial time.

\vspace{0.1in}

Another $t$-derivative gives 
$$\frac{\partial}{\partial t}(\frac{\partial^2 u}{\partial t^2})=
\Delta(\frac{\partial^2 u}{\partial t^2})+e^{-t}\<\tilde\omega
_t, \omega_0-\omega_\infty\>-\frac{\partial^2 u}{\partial 
t^2}-|\frac{\partial\widetilde\omega_t}{\partial t}|^2_{\tilde
\omega_t}.$$
Take summation with the first $t$-derivative to arrive at 
$$\frac{\partial}{\partial t}(\frac{\partial^2 u}{\partial t^2}+
\frac{\partial u}{\partial t})=\Delta(\frac{\partial^2 u}{\partial 
t^2}+\frac{\partial u}{\partial t})-(\frac{\partial^2 u}{\partial 
t^2}+\frac{\partial u}{\partial t})-|\frac{\partial\widetilde
\omega_t}{\partial t}|^2_{\widetilde\omega_t},$$
which gives
$$\frac{\partial^2 u}{\partial t^2}+\frac{\partial u}{\partial t}
\leqslant Ce^{-t}.$$
This implies the "essential decreasing" of volume form 
along the flow, i. e.    
$$\frac{\partial}{\partial t}(\frac{\partial u}{\partial t}+u)
\leqslant Ce^{-t},$$
which also induces
$$\frac{\partial u}{\partial t}\leqslant Ce^{-t}.$$

\vspace{0.1in}

Let's rewrite the metric flow equation as follows,
$${\rm Ric}(\widetilde\omega_t)=-\sqrt{-1}\partial\bar\partial
(u+\frac{\partial u}{\partial t})-\omega_\infty.$$

Taking trace with respect to $\widetilde\omega_t$ for the original 
metric flow equation and the one above, we have
$$R=e^{-t}\<\widetilde\omega_t, \omega_0-\omega_\infty\>-\Delta
(\frac{\partial u}{\partial t})-n=-\Delta(u+\frac{\partial u}{\partial 
t})-\<\widetilde\omega_t, \omega_\infty\>,$$
where $R$ denotes the scalar curvature of $\widetilde\omega_t$. 
Using the equations above, we also have
$$R=-n-\frac{\partial}{\partial t}(\frac{\partial u}{\partial t}+u),$$
and so the estimate got for $\frac{\partial}{\partial t}(\frac
{\partial u}{\partial t}+u)$ before is equivalent to the well 
known fact for scalar curvature. \\

We only consider smooth solution of K\"ahler-Ricci flow 
in $[0, T)\times X$ with finite time singularity at $T$. At 
this moment, we only need that the smooth limiting 
background form $\omega_T\geqslant 0$. It is essentially 
equivalent to assume $[\omega_T]$ has a smooth 
non-negative representative and presumably weaker 
than the class being "semi-ample", i.e. the existence of 
map $F$ before Theorem \ref{th:R-control}. 

\vspace{0.1in}

Recall the following equation used before 
$$\frac{\partial}{\partial t}\((1-e^{t-T})\frac{\partial u}
{\partial t}+u\)=\Delta\((1-e^{t-T})\frac{\partial u}
{\partial t}+u\)-n+\<\widetilde\omega_t, \omega_T\>$$
with the "$T$" in the equation chosen to the "$T$" 
above. With $\omega_T\geqslant 0$, by Maximum 
Principle, one has
$$(1-e^{t-T})\frac{\partial u}{\partial t}+u\geqslant -C.$$ 
Together with the upper bounds, we conclude 
$$|(1-e^{t-T})\frac{\partial u}{\partial t}+u|\leqslant C.$$

\vspace{0.1in}

\noindent $\bullet$ {\bf Parabolic Schwarz Estimate}

\vspace{0.1in}

Use the following setup as in \cite{song-tian} for the map 
$F$ before the statement of Theorem \ref{th:R-control}.  
Let $\varphi=\<\widetilde\omega_t, F^*\omega_{_{M}}\>$, 
then one has, over $[0, T)\times X$, 
$$(\frac{\partial}{\partial t}-\Delta)\varphi\leqslant \varphi
+C\varphi^2-H,$$
where  $C$ is related to the bisectional curvature bound 
of $\omega_{_{M}}$ near $F(X)$ and $H\geqslant 0$ is 
described as follows. Using normal coordinates locally 
over $X$ and $Y$, with indices $i, j$ and $\alpha, \beta$, 
$\varphi=|F_i^\alpha|^2$ and $H=|F_{ij}^\alpha|^2$ with 
all the summations. Notice that the normal coordinate 
over $X$ is changing along the flow with the metric.  
Using this inequality, one has 
$$(\frac{\partial}{\partial t}-\Delta){\rm log}\varphi\leqslant 
C\varphi+1.$$

\begin{remark}

For application, the map $F$ is coming from the class 
$[\omega_T]$ with $Y$ being some projective space 
$\mathbb{CP}^N$, and so $\omega_T$ is $F^*\omega$ 
where $\omega$ is (some mutiple of) Fubini-Study 
metric over $Y$. 

\end{remark}

Define 
$$v:=(1-e^{t-T})\frac{\partial u}{\partial t}+u$$ 
and we know $|v|\leqslant C$ for the previous step. 
We also have 
$$(\frac{\partial}{\partial t}-\Delta)v=-n+\<\widetilde\omega_t, 
\omega_T\>=-n+\varphi.$$
After taking a large enough positive constant $A$, 
the following inequality is true,
$$(\frac{\partial}{\partial t}-\Delta)({\rm log}\varphi-Av)
\leqslant -C\varphi+C.$$
Since $v$ is bounded, Maximum Principle can 
be used to deduce $\varphi\leqslant C$, i. e.  
$$\<\widetilde\omega_t, \omega_T\>\leqslant C.$$

\vspace{0.1in}

\noindent $\bullet$ {\bf Gradient and Laplacian Estimates}

\vspace{0.1in}

In this part, we derive gradient and Laplacian estimates 
for $v$. Recall that 
$$(\frac{\partial}{\partial t}-\Delta)v=-n+\varphi, ~~~ \varphi
=\<\widetilde\omega_t, \omega_T\>.$$
Standard computation gives:
$$(\frac{\partial}{\partial t}-\Delta)(|\nabla v|^2)=|\nabla v|^2-
|\nabla\nabla v|^2-|\nabla\bar\nabla v|^2+2{\rm Re}(\nabla
\varphi, \nabla v),$$
$$(\frac{\partial}{\partial t}-\Delta)(\Delta v)=\Delta v+({\rm 
Ric}(\widetilde\omega_t), \sqrt{-1}\partial\bar\partial v)+\Delta
\varphi.$$
Again, all the $\nabla$, $\Delta$ and $(\cdot, \cdot)$ are 
with respect to $\widetilde\omega_t$ and $\nabla\bar\nabla v$ 
is just $\partial\bar\partial v$. \\ 

Consider $\Psi:=\frac{|\nabla v|^2}{C-v}$. Since $v$ is 
bounded, one can easily make sure the denominator is 
positive, bounded and also away from $0$. We have the 
following computation,
\begin{equation}
\begin{split}
& ~~(\frac{\partial}{\partial t}-\Delta)\Psi= (\frac{\partial}
{\partial t}-\Delta)(\frac{|\nabla v|^2}{C-v}) \\
&= \frac{1}{C-v}\cdot\frac{\partial}{\partial t}(|\nabla v|^2)
+\frac{|\nabla v|^2}{(C-v)^2}\cdot\frac{\partial v}{\partial t}
-\( \frac{(|\nabla v|^2)_{\bar i}}{C-v}+\frac{v_{\bar i}
|\nabla v|^2}{(C-v)^2}\)_i \\
&= \frac{|\nabla v|^2}{(C-v)^2}\cdot(\frac{\partial}{\partial 
t}-\Delta)v+\frac{1}{C-v}\cdot(\frac{\partial}{\partial t}-\Delta)
(|\nabla v|^2)-\frac{v_i\cdot(|\nabla v|^2)_{\bar i}}{(C-v)^2}- 
v_{\bar i}\cdot\(\frac{|\nabla v|^2}{(C-v)^2}\)_i \\
&= \frac{|\nabla v|^2}{(C-v)^2}\cdot(\frac{\partial}{\partial 
t}-\Delta)v+\frac{1}{C-v}\cdot(\frac{\partial}{\partial t}-\Delta)
(|\nabla v|^2)-\frac{2{\rm Re}(\nabla v, \nabla |\nabla v|^2)}
{(C-v)^2}-\frac{2|\nabla v|^4}{(C-v)^3}. \nonumber
\end{split}
\end{equation} 
Plug in the results from before and rewrite the differential 
equality for $\Psi$ below,
\begin{equation}
\label{gradient2}
\begin{split}
& ~~(\frac{\partial}{\partial t}-\Delta)\Psi \\
&= \frac{(-n+\varphi)|\nabla v|^2}{(C-v)^2}+\frac{|\nabla v|^2-
|\nabla\nabla v|^2-|\nabla \bar\nabla v|^2}{C-v}+\frac{2{\rm 
Re}(\nabla\varphi, \nabla v)}{C-v} \\
& ~~~~-\frac{2{\rm Re}(\nabla v, \nabla |\nabla v|^2)}{(C-v)
^2}-\frac{2|\nabla v|^4}{(C-v)^3}. 
\end{split}
\end{equation}

The computations below are useful next.  
\begin{equation}
\begin{split}
|(\nabla v, \nabla |\nabla v|^2)|
&= |v_i(v_j v_{\bar j})_{\bar i}| \\
&= |v_iv_{\bar j}v_{j\bar i}+v_i
v_jv_{\bar j\bar i }| \\
&\leqslant |\nabla v|^2(|\nabla\nabla v|+|\nabla\bar\nabla 
v|) \\
&\leqslant \sqrt2 |\nabla v|^2(|\nabla\nabla v|^2+|\nabla
\bar\nabla v|^2)^{\frac{1}{2} }. \nonumber
\end{split}
\end{equation}
$$\nabla\Psi=\nabla\(\frac{|\nabla v|^2}{C-v}\)=\frac
{\nabla(|\nabla v|^2)}{C-v}+\frac{|\nabla v|^2\nabla 
v}{(C-v)^2}.$$

Together with the bounds for $\varphi$ and $C-v$, we can 
have the following computation with $\epsilon$ representing 
small positive constant (different from place to place),  
\begin{equation}
\begin{split}
& ~~(\frac{\partial}{\partial t}-\Delta)\Psi \\
&\leqslant C|\nabla v|^2+\epsilon\cdot|\nabla \varphi|^2-C
(|\nabla\nabla v|^2+|\nabla\bar\nabla v|^2)+ \\
& ~~~~~~ -(2-\epsilon){\rm Re}\(\nabla \Psi, \frac{\nabla v}
{C-v}\)-\epsilon\cdot\frac{{\rm Re}(\nabla v, \nabla |\nabla 
v|^2)}{(C-v)^2}-\epsilon\cdot\frac{|\nabla v|^4}{(C-v)^3} \\
&\leqslant C|\nabla v|^2+\epsilon\cdot|\nabla \varphi|^2-C
(|\nabla\nabla v|^2+|\nabla\bar\nabla v|^2)+ \\
& ~~~~~~ -(2-\epsilon){\rm Re}\(\nabla \Psi, \frac{\nabla v}
{C-v}\)+\epsilon\cdot(|\nabla\nabla v|^2+|\nabla\bar\nabla 
v|^2)-\epsilon\cdot|\nabla v|^4 \\ 
&\leqslant C|\nabla v|^2+\epsilon\cdot|\nabla \varphi|^2-
(2-\epsilon){\rm Re}\(\nabla \Psi, \frac{\nabla v}{C-v}\)-
\epsilon\cdot|\nabla v|^4. \nonumber
\end{split}
\end{equation}

We need a few more calculations to set up Maximum 
Principle argument. Recall that $\varphi=\<\widetilde\omega_t, 
\omega_T\>$ and, 
$$(\frac{\partial}{\partial t}-\Delta)\varphi\leqslant \varphi+
C\varphi^2-H.$$
With the description of $H$ before and the estimate for 
$\varphi$, i. e.  $\varphi\leqslant C$ from Schwarz estimate, 
we can conclude that 
$$H\geqslant C|\nabla\varphi|^2.$$ 
Now one arrives at
\begin{equation}
\label{schwarz2} 
(\frac{\partial}{\partial t}-\Delta)\varphi\leqslant C-C
|\nabla\varphi|^2.
\end{equation}
We also have
\begin{equation}
\label{elementary}
|\(\nabla \varphi, \frac{\nabla v}{C-v}\)|\leqslant \epsilon
\cdot|\nabla \varphi|^2+C\cdot|\nabla v|^2.
\end{equation}

\vspace{0.2in}

Now consider the function $\Psi+\varphi$. By choosing $
\epsilon>0$ small enough above (which also affects the 
choices of $C$'s), we have 
$$(\frac{\partial}{\partial t}-\Delta)(\Psi+\varphi)\leqslant C+
C|\nabla v|^2-\epsilon\cdot|\nabla v|^4-(2-\epsilon){\rm Re}
\(\nabla(\Psi+\varphi), \frac{\nabla v}{C-v}\).$$
At the maximum value point of $\Psi+\varphi$, we know 
$|\nabla v|^2$ can not be too large. It's then easy to 
conclude the upper bound for this term, and so for $\Psi$. 
Hence we have bounded the gradient, i. e.   
$$|\nabla v|\leqslant C.$$  

Now we want to do similar thing for the Laplacian, $\Delta 
v$. Define the funtion $\Phi:=\frac{C-\Delta v}{C-v}$. Similar computation as before gives the following
\begin{equation}
\label{laplacian2}
\begin{split}
& ~~ (\frac{\partial}{\partial t}-\Delta)\Phi= (\frac{\partial}
{\partial t}-\Delta)(\frac{C-\Delta v}{C-v}) \\
&= -\frac{1}{C-v}\cdot (\frac{\partial}{\partial t}-\Delta)
\Delta v+\frac{C-\Delta v}{(C-v)^2}\cdot(\frac{\partial}
{\partial t}-\Delta)v+\frac{2{\rm Re}(\nabla v, \nabla\Delta v)}
{(C-v)^2}+ \\
&~~~~~~ -\frac{2|\nabla v|^2(C-\Delta v)}{(C-v)^3} \\
&= -\frac{1}{C-v}\cdot(\Delta v+({\rm Ric}(\widetilde\omega_t), 
\sqrt{-1}\partial\bar\partial v)+\Delta \varphi)+\frac{C-
\Delta v}{C-v}\cdot(-n+\varphi) \\
&~~~~~~ +\frac{2{\rm Re}(\nabla v, \nabla\Delta v)}
{(C-v)^2}-\frac{2|\nabla v|^2(C-\Delta v)}{(C-v)^3}. 
\end{split}
\end{equation}

We also have $\nabla(\frac{C-\Delta v}{C-v})=\frac
{(C-\Delta v)\nabla v}{(C-v)^2}-\frac{\nabla\Delta v}
{C-v}$. Recall that it is already known $(0\leqslant) 
\varphi\leqslant C$. The following inequality follows 
from standard computation as in \cite{song-tian} 
and has actually been used for parabolic Schwarz 
estimate, 
$$\Delta \varphi\geqslant ({\rm Ric}(\widetilde\omega_t), 
\omega_T)+H-C\varphi^2,$$
where $H\geqslant C|\nabla \varphi|^2\geqslant 0$ 
from the bound of $\varphi$ as mentioned before. 
Now we have 
\begin{equation}
\label{cases}
({\rm Ric}(\widetilde\omega_t), \sqrt{-1}\partial\bar\partial 
v)+\Delta \varphi\geqslant ({\rm Ric}(\widetilde\omega_t), 
\sqrt{-1}\partial\bar\partial v+\omega_T)-C.
\end{equation}

\vspace{0.1in}

We are considering the case $T<\infty$. Recall that 
$v=(1-e^{t-T})\frac{\partial u}{\partial t}+u$ and $
\omega_T=\omega_\infty+e^{-T}(\omega_0-\omega
_\infty)$. We have  
\begin{equation*}
\begin{split}
{\rm Ric}(\widetilde\omega_t)
&= -\sqrt{-1}\partial\bar\partial \(\frac{\p u}{\p t}+u\)-
\omega_\infty \\
&= -\sqrt{-1}\p\bar\p v -\omega_T-e^{t-T}\sqrt{-1}\p\bar
\p \frac{\p u}{\p t}+e^{-T}(\omega_0-\omega_\infty) \\
&=  -\sqrt{-1}\p\bar\p v -\omega_T-e^{t-T}\(\sqrt{-1}
\p\bar\p \frac{\p u}{\p t}-e^{-t}(\omega_0-\omega_\infty)
\) \\
&=  -\sqrt{-1}\p\bar\p v -\omega_T-e^{t-T}\frac{\p \tilde
\omega_t}{\p t} \\
&=  -\sqrt{-1}\p\bar\p v -\omega_T-e^{t-T}\(-\Ric(\tilde
\omega_t)-\widetilde\omega_t\),  
\end{split}
\end{equation*}
which gives 
$$(1-e^{t-T})\Ric(\widetilde\omega_t)=-\sqrt{-1}\p\bar\p v-
\omega_T+e^{t-T}\widetilde\omega_t,$$
and so 
$$\Ric(\widetilde\omega_t)=-\frac{\sqrt{-1}\p\bar\p v+\omega
_T}{1-e^{t-T}}+\frac{e^{t-T}}{1-e^{t-T}}\widetilde\omega_t,$$
$$(1-e^{t-T})R=-\Delta v-\<\widetilde\omega_t, \omega_T\>
+ne^{t-T}.$$

\vspace{0.1in}

As $R\geqslant -C$ and $0\leqslant \<\widetilde\omega_t, 
\omega_T\>$, we have $\Delta v\leqslant C$. \\ 

Now we can continue the estimation (\ref{cases}) as 
follows.   
\begin{equation*}
\begin{split}
&~~ ({\rm Ric}(\widetilde\omega_t), \sqrt{-1}\partial\bar\partial 
v)+\Delta \varphi \\
&\geqslant ({\rm Ric}(\widetilde\omega_t), \sqrt{-1}\partial\bar
\partial v+\omega_T)-C \\
&\geqslant \(-\frac{\sqrt{-1}\p\bar\p v+\omega_T}{1-e^{t-T}}
+\frac{e^{t-T}}{1-e^{t-T}}\widetilde\omega_t,  \sqrt{-1}\partial
\bar\partial v+\omega_T\)-C \\
&= -\frac{|\sqrt{-1}\p\bar\p v+\omega_T|^2}{1-e^{t-T}}+
\frac{e^{t-T}(\Delta v+\<\widetilde\omega_t, \omega_T\>)}
{1-e^{t-T}}-C. 
\end{split}
\end{equation*}
As $C^{-1}(T-t)\leqslant 1-e^{t-T}\leqslant C(T-t)$ for 
$t\in [0, T)$, using $\Delta v=\Delta v-C+C$ and $0
\leqslant\<\widetilde\omega_t, \omega_T\>\leqslant C$, 
we have 
\begin{equation*}
\begin{split}
&~~ ({\rm Ric}(\widetilde\omega_t), \sqrt{-1}\partial\bar
\partial v)+\Delta \varphi \\
&= -\frac{|\sqrt{-1}\p\bar\p v+\omega_T|^2}{1-e^{t-T}}+
\frac{e^{t-T}(\Delta v+\<\widetilde\omega_t, \omega_T\>)}
{1-e^{t-T}}-C \\ 
&\geqslant -\frac{1}{T-t}\((1+\epsilon)|\sqrt{-1}\p\bar\p 
v|^2+C|\omega_T|^2\)-\frac{C}{T-t}(C-\Delta v)-\frac{C}
{T-t} \\
&\geqslant -\frac{1+\epsilon}{T-t}|\sqrt{-1}\p\bar\p v|^2-
\frac{C}{T-t}(C-\Delta v)-\frac{C}{T-t}. 
\end{split}
\end{equation*}
  
Now we can continue the computation for $\Phi$, 
(\ref{laplacian2}) as follows.  
\begin{equation*}
(\frac{\partial}{\partial t}-\Delta)\Phi\leqslant \frac{C}{T-t}+
\frac{C}{T-t}\cdot(C-\Delta v)+\frac{(1+\epsilon)|\nabla\bar
\nabla v|^2}{(T-t)(C-v)}-2{\rm Re}\(\nabla \Phi, \frac{\nabla 
v}{C-v}\). 
\end{equation*}

Using $\Phi=\frac{C-\Delta v}{C-v}\geqslant C(C-\Delta v)$, 
one arrives at   
\begin{equation*}
(\frac{\partial}{\partial t}-\Delta)((T-t)\Phi)\leqslant C+C\cdot 
(C-\Delta v)+\frac{(1+\epsilon)|\nabla\bar\nabla v|^2}{C-v}-
2{\rm Re}(\nabla \((T-t)\Phi\), \frac{\nabla v}{C-v})  
\end{equation*}
In sight of (\ref{schwarz2}) and (\ref{elementary}), we have,  
\begin{equation*}
(\frac{\partial}{\partial t}-\Delta)\varphi\leqslant C-4{\rm Re}
\(\nabla\varphi, \frac{\nabla v}{C-v}\)+C|\nabla v|^2. 
\end{equation*}
Also, (\ref{gradient2}) can be rewritten as 
\begin{equation}
\begin{split}
(\frac{\partial}{\partial t}-\Delta)\Psi 
&\leqslant \frac{(-n+\varphi)|\nabla v|^2}{(C-v)^2}+\frac
{|\nabla v|^2-|\nabla\nabla v|^2-|\nabla \bar\nabla v|^2}
{C-v}+ \\
&~~~~~~ 2{\rm Re}\(\nabla\varphi, \frac{\nabla v}{C-v}\)
-2{\rm Re}\(\Psi, \frac{\nabla v}{C-v}\). 
\nonumber
\end{split}
\end{equation}

Using the bound for $|\nabla v|$ and choosing $\epsilon<1$, 
we have 
\begin{equation*}
\begin{split}
&~~ (\frac{\partial}{\partial t}-\Delta)\((T-t)\Phi+2\Psi+2
\varphi\) \\ 
&\leqslant C+C\cdot(C-\Delta v)-2{\rm Re}\(\nabla\((T-t)
\Phi+2\Psi+2\varphi\), \frac{\nabla v}{C-v}\)-C|\nabla\bar
\nabla v|^2 \\
&\leqslant  C+C\cdot(C-\Delta v)-2{\rm Re}\(\nabla\((T-t)
\Phi+2\Psi+2\varphi\), \frac{\nabla v}{C-v}\)-C(C-\Delta 
v)^2
\end{split}
\end{equation*}
where $|\nabla\bar\nabla v|^2\geqslant C(\Delta v)^2
\geqslant C(C-\Delta v)^2-C$ is used for the second 
$\leqslant$. \\ 

Now we apply Maximum Principle. At maximum value 
point of the function $(T-t)\Phi+2\Psi+2\varphi$, we have 
$C-\Delta v\leqslant C$. Using the bounds on $\Psi$ and 
$\varphi$, $(T-t)\Phi+2\Psi+2\varphi\leqslant C$ and so 
$$\Phi\leqslant \frac{C}{T-t}, ~~\text{i.e.}~~ \Delta v
\geqslant -\frac{C}{T-t}.$$
Finally since $(1-e^{t-T})R=-\Delta v-\<\widetilde\omega_t, 
\omega_T\>+ne^{t-T}$, we conclude that 
$$R\leqslant \frac{C}{(T-t)^2}.$$

\vspace{0.1in}

Theorem \ref{th:R-control} is proved. 

\section{Further Remarks}

There are several closely related results worth mentioning. 
The last two remarks below should give people the idea 
about the essential difference between finite time and infinite 
time singular case for K\"ahler-Ricci flow. 

\begin{itemize}

\item In \cite{sesum-tian}, following Perelman's idea, Sesum 
and Tian proved that for $X$ with $c_1(X)>0$, for any 
initial K\"ahler metric $\omega$ such that $[\omega]=
c_1(X)$, the K\"ahler-Ricci flow 
$$\frac{\p \widehat\omega_t}{\p t}=-\Ric(\widehat\omega
_t)+\widehat\omega_t$$
has uniformly bounded scalar curvature and diameter 
for $\widehat\omega_t$ where $t\in [0, \infty)$. Using 
simply rescaling of time and metric, one can see for 
our flow (\ref{krf}) with $[\omega_0]=c_1(X)$, 
$$R(\widetilde\omega_t)\leqslant \frac{C}{T-t}$$ 
for $t\in [0, T)$ where the finite singular time $T=\log 2$, 
which is a better control than Theorem \ref{th:R-control} 
for this special case. 

\item For the infinite time limiting case, the scalar curvature 
would be bounded if the infinite time limiting class provides 
a smooth holomorphic fibration for $X$, i.e. the map $F$ as 
in our setting is a smooth fibration. This is actually proved in 
\cite{song-tian} if one only considers smooth collapsing case. 

\item For the infinite time limiting case, the scalar curvature 
would also be bounded if the limiting class is "semi-ample 
and big", i.e. the (possibly singular) image of the map $F$ 
is of the same dimension as $X$. This result is proved in 
\cite{b-curvature}. The more recent work of Yuguang Zhang, 
\cite{miyaoka-yau}, has given a nice application of it.   

\end{itemize}

\vspace{0.4in}

Zhou Zhang 

Department of Mathematics

University of Michigan, at Ann Arbor 

MI 48109

zhangou@umich.edu


\begin{thebibliography}{$$}

\bibitem{cascini-lanave} Cascini, Paolo; La Nave, Gabriele: 
K\"ahler-Ricci flow and the minimal model program for 
projective varieties. ArXiv:math/0603064 (math.AG).

\bibitem{demailly-paun} Demailly, Jean-Pierre; Paun,
Mihai: Numerical characterization of the K\"ahler cone
of a compact K\"ahler manifold. Ann. of Math. 159
(2004), 1247--1274.

\bibitem{ham} Hamilton, Richard S.: Three-manifolds
with positive Ricci curvature. J. Differential Geom. 17
(1982), no. 2, 255-306.

\bibitem{ham4} Hamilton, Richard S.: The formation of
singularities in the Ricci flow.  Surveys in differential 
geometry, Vol. II (Cambridge, MA, 1993),  7--136, Int. 
Press, Cambridge, MA, 1995.

\bibitem{sesum} Sesum, Natasa: Curvature tensor 
under the Ricci flow.  Amer. J. Math.  127  (2005),  
no. 6, 1315--1324.

\bibitem{sesum-tian} Sesum, Natasa; Tian, Gang: 
Bounding scalar curvature and diamater along the 
K\"ahler Ricci flow (after Perelman).  J. Inst. Math. 
Jussieu  7  (2008), no. 3, 575--587. 53C44.

\bibitem{song-tian} Song, Jian; Tian, Gang: The 
K\"ahler-Ricci flow on surfaces of positive Kodaira 
dimension. Invent. Math. 170 (2007), no. 3, 609--653.  

\bibitem{t-znote} Tian, Gang; Zhang, Zhou: On the 
K\"ahler-Ricci flow on projective manifolds of general 
type. Chinese Annals of Mathematics - Series B, 
Volume 27, Number 2, 179--192.   

\bibitem{miyaoka-yau} Zhang, Yuguang: Miyaoka-Yau 
inequality for minimal projective manifolds of general 
type. ArXiv:0812.0462 (math.DG) (math.AG). 

\bibitem{b-curvature} Zhang, Zhou: Scalar Curvature 
Bound for K\"ahler-Ricci Flows over Minimal Manifolds 
of General Type. ArXiv:0801.3248 (math.DG).  

\end{thebibliography}
\end{document}